\documentclass[a4paper,12pt]{amsart}
\usepackage{amsfonts}
\usepackage{amssymb}
\usepackage{ifthen}
\usepackage{graphicx}
\usepackage{color}
\nonstopmode \numberwithin{equation}{section}
\newtheorem{thm}{Theorem}
\newtheorem{lem}{Lemma}
\newtheorem{cor}{Corollary}

\newtheorem{cl}{Claim}
\newtheorem{ca}{Case}
\newtheorem{sca}{Subcase}
\newtheorem{scl}{Subclaim}
\newtheorem{conj}[equation]{Conjecture}

\theoremstyle{definition}
\newtheorem{defn}{Definition}

\newtheorem{op}[equation]{Open Problem}
\newtheorem{ques}[equation]{Question}
\newtheorem{rem}{Remark}[section]
\newtheorem{exam}[equation]{Example}

\newcounter {own}
\def\theown {\thesection       .\arabic{own}}

\newenvironment{pf}[1][]{%
 \vskip 3mm
 \noindent
 \ifthenelse{\equal{#1}{}}%
  {{\slshape Proof. }}%
  {{\slshape #1.} }%
 }%
{\qed\bigskip}

\newcounter{alphabet}
\newcounter{tmp}
\newenvironment{Thm}[1][]{\refstepcounter{alphabet}%
\bigskip%
\noindent%
{\bf Theorem \Alph{alphabet}}%
\ifthenelse{\equal{#1}{}}{}{ (#1)}%
{\bf .} \itshape}{\vskip 8pt}

\makeatletter
\newcommand{\Ref}[1]{\@ifundefined{r@#1}{}{\setcounter{tmp}{\ref{#1}}\Alph{tmp}}}
\makeatother

\newenvironment{Lem}[1][]{\refstepcounter{alphabet}%
\bigskip%
\noindent%
{\bf Lemma \Alph{alphabet}}%
{\bf .} \itshape}{\vskip 8pt}

\newcommand{\ID}{{\mathbb D}}

\def\be{\begin{equation}}
\def\ee{\end{equation}}

\newcommand{\bee}{\begin{enumerate}}
\newcommand{\eee}{\end{enumerate}}

\newcommand{\blem}{\begin{lem}}
\newcommand{\elem}{\end{lem}}
\newcommand{\bthm}{\begin{thm}}
\newcommand{\ethm}{\end{thm}}
\newcommand{\bcor}{\begin{cor}}
\newcommand{\ecor}{\end{cor}}
\newcommand{\beg}{\begin{exam}}
\newcommand{\eeg}{\end{exam}}
\newcommand{\begs}{\begin{examples}}
\newcommand{\eegs}{\end{examples}}
\newcommand{\bdefe}{\begin{defn}}
\newcommand{\edefe}{\end{defn}}
\newcommand{\bprob}{\begin{prob}}
\newcommand{\eprob}{\end{prob}}
\newcommand{\bques}{\begin{ques}}
\newcommand{\eques}{\end{ques}}
\newcommand{\bei}{\begin{itemize}}
\newcommand{\eei}{\end{itemize}}
\newcommand{\bcon}{\begin{conj}}
\newcommand{\econ}{\end{conj}}
\newcommand{\bop}{\begin{op}}
\newcommand{\eop}{\end{op}}

\newcommand{\bca}{\begin{ca}}
\newcommand{\eca}{\end{ca}}
\newcommand{\bsca}{\begin{sca}}
\newcommand{\esca}{\end{sca}}

\newcommand{\bcl}{\begin{cl}}
\newcommand{\ecl}{\end{cl}}

\newcommand{\bscl}{\begin{scl}}
\newcommand{\escl}{\end{scl}}

\newcommand{\bcons}{\begin{conjs}}
\newcommand{\econs}{\end{conjs}}
\newcommand{\bprop}{\begin{propo}}
\newcommand{\eprop}{\end{propo}}
\newcommand{\br}{\begin{rem}}
\newcommand{\er}{\end{rem}}
\newcommand{\brs}{\begin{rems}}
\newcommand{\ers}{\end{rems}}
\newcommand{\bo}{\begin{obser}}
\newcommand{\eo}{\end{obser}}
\newcommand{\bos}{\begin{obsers}}
\newcommand{\eos}{\end{obsers}}
\newcommand{\bpf}{\begin{pf}}
\newcommand{\epf}{\end{pf}}
\newcommand{\ba}{\begin{array}}
\newcommand{\ea}{\end{array}}
\newcommand{\beq}{\begin{eqnarray}}
\newcommand{\beqq}{\begin{eqnarray*}}
\newcommand{\eeq}{\end{eqnarray}}
\newcommand{\eeqq}{\end{eqnarray*}}

\newcommand{\ds}{\displaystyle}

\newcounter{minutes}
\divide\time by 60
\newcounter{hours}
\multiply\time by 60 \addtocounter{minutes}{-\time}

\begin{document}

\bibliographystyle{amsplain}
\title [] {
 Lipschitz spaces and bounded mean oscillation of
harmonic mappings}

\def\thefootnote{}
\footnotetext{ \texttt{\tiny File:~\jobname .tex,
          printed: \number\day-\number\month-\number\year,
          \thehours.\ifnum\theminutes<10{0}\fi\theminutes}
} \makeatletter\def\thefootnote{\@arabic\c@footnote}\makeatother


\author{Sh. Chen}
\address{Sh. Chen, Department of Mathematics and Computational
Science, Hengyang Normal University, Hengyang, Hunan 421008,
People's Republic of China.} \email{shlchen1982@yahoo.com.cn}


\author{S.  Ponnusamy $^\dagger $
}
\address{S. Ponnusamy,
Indian Statistical Institute (ISI), Chennai Centre, SETS (Society
for Electronic Transactions and security), MGR Knowledge City, CIT
Campus, Taramani, Chennai 600 113, India. }
\email{samy@isichennai.res.in, samy@iitm.ac.in}

\author{M. Vuorinen}
\address{M. Vuorinen, Department of Mathematics,
University of Turku,  Turku 20014, Finland.} \email{vuorinen@utu.fi}

\author{X. Wang${}^{~\mathbf{*}}$}
\address{X. Wang, Department of Mathematics,
Hunan Normal University, Changsha, Hunan 410081, People's Republic
of China.} \email{xtwang@hunnu.edu.cn}

\keywords{Harmonic  mapping, majorant, Lipschitz space, $BMO_{p}$,
equivalent
modulus, harmonic Bloch space, Green's theorem.\\
$
^\dagger$ {\tt This author is on leave from the Department of Mathematics,
Indian Institute of Technology Madras, Chennai-600 036, India}
}

\begin{abstract}
In this paper, we first study the bounded mean oscillation of planar
harmonic mappings, then a relationship between  Lipschitz-type spaces and
equivalent modulus of real harmonic mappings is established. At last, we obtain sharp estimates on Lipschitz number of planar harmonic
mappings in terms of
bounded mean oscillation norm, which shows that the harmonic Bloch space is
isomorphic to $BMO_{2}$ as a Banach space.
\end{abstract}


\maketitle \pagestyle{myheadings} \markboth{ SH. Chen, S. Ponnusamy,
M. Vuorinen and X. Wang }{ Lipschitz spaces and bounded mean
oscillation of harmonic mappings }

\section{Introduction and main results}\label{csw-sec1}
Let $\mathbb{C}$ denote the complex plane. For $a\in\mathbb{C}$, let
$\ID(a,r)=\{z:\, |z-a|<r\}$. In particular, we use $\mathbb{D}_r$ to
denote the disk $\mathbb{D}(0,r)$ and  $\mathbb{D}$ the unit disk $\ID_1$.
A complex-valued function $f$ defined on $\ID$ is called {\it
harmonic} in $\ID$ if and only if both the real and the imaginary parts of $f$
are real harmonic in $\ID$.
It is known that every harmonic mapping $f$ defined in $\ID$ admits
a  decomposition $f=h+\overline{g}$, where $h$ and $g$ are analytic
in $\ID$. We refer to \cite{Clunie-Small-84,Co, Du,Ha,S} for the
theory of planar harmonic mappings.
For harmonic mappings $f$ defined on $\mathbb{D}$, we use the following
standard notations:
$$\Lambda_{f}(z)=\max_{0\leq \theta\leq 2\pi}|f_{z}(z)+e^{-2i\theta}f_{\overline{z}}(z)|
=|f_{z}(z)|+|f_{\overline{z}}(z)|
$$
and
$$\lambda_{f}(z)=\min_{0\leq \theta\leq 2\pi}|f_{z}(z)+e^{-2i\theta}f_{\overline{z}}(z)|
=\big | \, |f_{z}(z)|-|f_{\overline{z}}(z)|\, \big |.
$$

A continuous increasing function $\omega:\, [0,\infty)\rightarrow
[0,\infty)$ with $\omega(0)=0$ is called a {\it majorant} if
$\omega(t)/t$ is non-increasing for $t>0$ (see \cite{D,P}). Given a
subset $\Omega$ of $\mathbb{C}$, a function $f:\, \Omega\rightarrow
\mathbb{C}$ is said to belong to the {\it Lipschitz space
$L_{\omega}(\Omega)$} if there is a positive constant $M$ such that
\be\label{eq1x}
|f(z)-f(w)|\leq M\omega(|z-w|) ~\mbox{ for all $z,\ w\in\Omega$.}
\ee

For $\delta_{0}>0$ and $0<\delta<\delta_{0}$, we consider the following conditions on
a majorant $\omega$:
\be\label{eq2x}
\int_{0}^{\delta}\frac{\omega(t)}{t}\,dt\leq M \omega(\delta)
\ee
and
\be\label{eq3x}
\delta\int_{\delta}^{+\infty}\frac{\omega(t)}{t^{2}}\,dt\leq M \omega(\delta),
\ee
where $M$ denotes a positive constant.

A majorant $\omega$ is said to be {\it regular} if it satisfies the
conditions (\ref{eq2x}) and (\ref{eq3x}) (see \cite{D,P}). 


Dyakonov \cite{D} discussed the relationship between the Lipschitz
space and the bounded mean oscillation on holomorphic functions in
$\mathbb{D}$, and obtained the following result.
 In order to state an analogue of
Theorem \Ref{ThmA} for planar harmonic mappings, we first introduce
some notation. Let $G$ be a domain of $\mathbb{C}.$  We use
$d_{G}(z)$ to denote the Euclidean distance from $z$ to the boundary
$\partial G$ of $G$. In particular, we always use $d(z)$ to denote
the Euclidean distance from $z$ to the boundary of $\mathbb{D}$.

\begin{Thm}{\rm \cite[Theorem 1]{D}}\label{ThmA}
Suppose that $f$ is a holomorphic function in $\mathbb{D}$ which is
continuous up to the boundary of $\mathbb{D}$. If $\omega$ and
$\omega^{2}$ are regular majorants,  then
$$f\in L_{\omega}(\mathbb{D})\Longleftrightarrow P_{|f|^{2}}(z)-|f(z)|^{2}\leq M\omega^{2}(d(z)),
$$
where
$$P_{|f|^{2}}(z)=\frac{1}{2\pi}\int_{0}^{2\pi}\frac{1-|z|^{2}}{|z-e^{i\theta}|^{2}}|f(e^{i\theta})|^{2}\,d\theta.
$$
\end{Thm}







 The following result is an analogue
of Theorem \Ref{ThmA} for planar harmonic mappings.

\begin{thm}\label{thm1}
Suppose that $\omega$ is a majorant and that $f$ is a harmonic mapping in
$\mathbb{D}$. Then $\Lambda_{f}(z)\leq
M\omega\Big(\frac{1}{d(z)}\Big)$ in $\mathbb{D}$ if and only if for every
$r\in(0,1-|z|]$,
$$\frac{1}{|\mathbb{D}(z,r)|}\int_{\mathbb{D}(z,r)}|f(\zeta)-f(z)|\,dA(\zeta)\leq Mr\omega\big(\frac{1}{r}\big),
$$
where 
$dA$ denotes the area
measure in $\mathbb{D}$.
\end{thm}

\begin{defn}\label{defn2}
Let $f$ be harmonic in $\mathbb{D}$. For $p\in[1,\infty)$, we say
$f\in BMO_{p}$ if
$$\|f\|_{BMO_{p}}=\sup_{\mathbb{D}(z,r)\subseteq\mathbb{D}}
\left\{\frac{1}{|\mathbb{D}(z,r)|}\int_{\mathbb{D}(z,r)}\left|f(\zeta)-
\frac{1}{|\mathbb{D}(z,r)|}\int_{\mathbb{D}(z,r)}f(\xi)\,dA(\xi)\right|^{p}\,dA(\zeta)\right\}^{1/p}
$$
is bounded, where   $r\in(0,1-|z|]$.
\end{defn}


In particular, by taking $\omega(t)=t$ in Theorem \ref{thm1}, we get the following
result.

\begin{cor}\label{cor-1}
Let $f$ be a harmonic mapping in $\mathbb{D}$. Then $f\in BMO_{1}$ if
and only if $\Lambda_{f}(z)\leq M\frac{1}{d(z)}$ holds in $\mathbb{D}$.
\end{cor}

In \cite{D}, Dyakonov also investigated the property of equivalent
modulus for holomorphic functions in $\mathbb{D}$ and obtained

\begin{Thm}{\rm \cite[Theorem 2]{D}}\label{ThmB}
Let $\omega$ be a regular majorant and $f$ be a holomorphic function
in $\mathbb{D}$ and continuous up to the boundary $\partial
\mathbb{D}$. Then
$$f\in L_{\omega}(\mathbb{D})\Longleftrightarrow |f|\in L_{\omega}(\mathbb{D})\Longleftrightarrow
|f|\in L_{\omega}(\mathbb{D},\partial \mathbb{D}),
$$ where $L_{\omega}(\mathbb{D},\partial \mathbb{D})$ denotes the class of continuous
functions $F$ on $\mathbb{D}\cup\partial \mathbb{D}$ which satisfy
{\rm (\ref{eq1x})} with some positive constant $C$, whenever $z\in
\mathbb{D}$ and $w\in\partial \mathbb{D}$.
\end{Thm}

 Later in \cite[Theorems A]{P},
Pavlovi${\rm \acute{c}}$ came up with a relatively simple proof of
the results of Dyakonov. Recently, many authors considered this
topic and generalized Dyakonov's results to quasiconformal mappings
and real harmonic functions in several variables for some special
majorant $\omega(t)= t^{\alpha}$, where $\alpha>0$ (see
\cite{ABM,D1,KP,M3,MM,Pav1,Pav2,Pav3}). For the general majorant
$\omega$ to holomorphic mappings and pluriharmonic mappings in the
unit ball, see \cite{CPW7,D1,Q}.

We will prove  the analog of Theorem \Ref{ThmB} for real harmonic
functions in the following form.



\begin{thm}\label{Thm-B}
Suppose that $\omega$ is a  majorant satisfying {\rm (\ref{eq2x})},
and that $G$ is a $L_{\omega}$-extension domain. If $f$ is a real
harmonic function in $G$ and continuous up to the boundary $\partial
G$, then
$$f\in L_{\omega}(G)\Longleftrightarrow |f|\in L_{\omega}(G)\Longleftrightarrow |f|\in L_{\omega}(G,\partial G),
$$ 
where $L_{\omega}(G,\partial G)$ denotes the class of continuous
functions $F$ on $G\cup\partial G$ which satisfy {\rm (\ref{eq1x})}
with some positive constant $C$, whenever $z\in G$ and $w\in\partial
G$.
\end{thm}

 Here a proper subdomain $G$ of $\mathbb{C}$ or
$\mathbb{R}^{2}$  is said to be {\it $L_{\omega}$-extension} if
$L_{\omega}(G)=\mbox{loc}L_{\omega}(G)$, where
$\mbox{loc}L_{\omega}(G)$ denotes the set of all functions $f:\,
G\rightarrow \mathbb{C}$ satisfying  (\ref{eq1x}) with a fixed
positive constant $M$, whenever $z\in G$ and $w\in G$ such that
$|z-w|<\frac{1}{2}d_{G}(z)$.  Obviously, the unit disk $\mathbb{D}$
is a $L_{\omega}$-extension domain.

In \cite{L}, the author proved that $G$ is a $L_{\omega}$-extension
domain if and only if each pair of points $z,w\in G$ can be joined
by a rectifiable curve $\gamma\subset G$ satisfying
\be\label{eq1.0}
\int_{\gamma}\frac{\omega(d_{G}(z))}{d_{G}(z)}\,ds(z) \leq
M\omega(|z-w|)
\ee
with some fixed positive constant $M=M(G,\omega)$, where $ds$ stands for the arc
length measure on $\gamma$.
See \cite{GM,L} for more details on $L_{\omega}$-extension domains.

We remark that in Theorem \ref{Thm-B}, we replace ``the unit disk
$\mathbb{D}$" and ``the regular majorant" in Theorem \Ref{ThmB} by
``a $L_{\omega}$-extension domain" and ``a majorant satisfying
\eqref{eq2x}, but not necessarily \eqref{eq3x}", respectively. In
fact, by using   \cite[Lemma A, Theorem 4, Corollary 2]{Pav2} and the similar proof method of Theorem
\ref{Thm-B}, we can prove that Theorem \ref{Thm-B} also holds for
real harmonic functions in the unit ball $\mathbb{B}^{n}$ of
$\mathbb{R}^{n}$.

For planar harmonic mappings, we obtain the following result which
is a generalization of  Theorem \Ref{ThmB}.

\begin{thm}\label{Thmy1}
Let $\omega$ be a  majorant satisfying {\rm (\ref{eq2x})} and $G$ be
a $L_{\omega}$-extension domain. Let $f=h+\overline{g}$ be a
harmonic mapping in $G$, where $g$ and $h$ are analytic functions in
$G$. Then $$f\in L_{\omega}(G)\Longleftrightarrow g,\ h\in
L_{\omega}(G)\Longleftrightarrow |g|,\ |h|\in L_{\omega}(G).$$
\end{thm}

\begin{defn}\label{defn-2}
A planar harmonic mapping $f$ in $\mathbb{D}$ is called a {\it
harmonic Bloch mapping} if
$$\beta_{f}=\sup_{ z,w\in\mathbb{D},\ z\neq w}\frac{|f(z)-f(w)|}{\rho(z,w)}<\infty.
$$
Here $\beta_{f}$ is called the {\it  Lipschitz number} of $f$ and
$$\rho(z,w)=\frac{1}{2}\log\left(\frac{1+|\frac{z-w}{1-\overline{z}w}|}
{1-|\frac{z-w}{1-\overline{z}w}|}\right)=\mbox{arctanh}\Big|\frac{z-w}{1-\overline{z}w}\Big|
$$
denotes the hyperbolic distance between $z$ and $w$ in $\mathbb{D}$.
\end{defn}

It is known that
$$\beta_{f}=\sup_{z\in\mathbb{D}}\big\{(1-|z|^{2})\Lambda_{f}(z)\big\}.
$$
Clearly, a harmonic Bloch mapping $f$ is uniformly continuous as
a map between metric spaces
$$
f: (\mathbb{D}, \rho) \to ( \mathbb{C}, | \cdot|)
$$
and for all $z,w \in \mathbb{D}$ we have the Lipschitz inequality
$$
|f(z)-f(w)| \le \beta_{f} \, \rho(z,w) \,.
$$
The reader is referred to \cite[Theorem 2]{Co} (or
\cite{CPW0,CPW1,CPW2}) for a proof. Then the set of all harmonic
Bloch mappings in $\mathbb{D}$ forms a harmonic Bloch space which is
denoted by $\mathcal{B}_{h}$. Uniform continuity with respect to a hyperbolic
metric is a central theme in \cite{V1,V2}

In \cite{CM,HP,Po}, the authors provided several characterizations
of $BMO_{2}$ on holomorphic functions. For the extensive discussions
on $BMO_{2}$, see \cite{CR,FS,HY,J,KO}.
In this paper, we will use $BMO_{2}$ norm to obtain a sharp estimate
on  harmonic Bloch mappings, which shows that $\mathcal{B}_{h}$ is
isomorphic to $BMO_{2}$ as a Banach space. Our result is given
below.

\begin{thm}\label{thm-2}
If $f$ is harmonic in $\mathbb{D}$, then
\be\label{eq-a}
\|f\|_{BMO_{2}}\leq\beta_{f}\leq 2\|f\|_{BMO_{2}}.
\ee
Moreover, the estimates of \eqref{eq-a} are sharp. The extreme
harmonic mappings of the first inequality are constant functions, and the
extreme harmonic mappings of the second inequality are the mappings with the form
$f(z)=C(z+\overline{z})$, where $C$ denotes a constant.
\end{thm}

The proofs of Theorems \ref{thm1} and \ref{Thm-B} will be presented
in Section \ref{csw-sec2}, and the proof of Theorem \ref{thm-2} will
be given in Section \ref{csw-sec3}.

\section{ Bounded mean oscillation and equivalent modulus }\label{csw-sec2}
The following lemma easily follows from a simple computation (cf. \cite{CPW6}).

\begin{lem}
\label{Lem1}
Let $f$ be a complex-valued continuously differentiable function
defined on $\mathbb{D}$ and $f=u+iv$, where $u$ and $v$ are
real-valued functions. Then for $z=x+iy\in\mathbb{D}$,
\be\label{eqs1}
\Lambda_{f}(z)\leq |\nabla u(x,y)|+|\nabla v(x,y)|,
\ee
where $\nabla u=(u_{x},u_{y})$ and $\nabla v=(v_{x},v_{y})$.
\end{lem}

 Then we have

\begin{lem}\label{lem1}
Suppose $f$ is a harmonic mapping in
$\overline{\mathbb{D}}(a,r)$, where $r$ is a positive constant. Then
$$\Lambda_{f}(a)\leq\frac{2}{\pi r}\int_{0}^{2\pi}|f(a)-f(a+re^{i\theta})|\,d\theta.
$$
\end{lem}
\bpf Let $f=u+iv$ be a harmonic mapping in
$\overline{\mathbb{D}}(a,r)$, where $u$ and $v$ are real harmonic
functions. Without loss of generality, we may assume that $a=0$ and
$f(0)=0.$ By Poisson's formula, we have
$$u(z)=\frac{1}{2\pi}\int_{0}^{2\pi}\frac{r^{2}-|z|^{2}}{|z-re^{i\theta}|^{2}}u(re^{i\theta})\,d\theta, \quad |z|<r.
$$
By calculations, we get ($z=x=iy$)
$$u_{x}(z)=\frac{1}{2\pi}\int_{0}^{2\pi}\frac{-2x|z-re^{i\theta}|^{2}-2(r^{2}-|z|^{2})(x-r\cos\theta)}
{|z-re^{i\theta}|^{4}}u(re^{i\theta})\,d\theta
$$
and similarly
$$u_{y}(z)=\frac{1}{2\pi}\int_{0}^{2\pi}\frac{-2y|z-re^{i\theta}|^{2}-2(r^{2}-|z|^{2})(y-r\sin\theta)}
{|z-re^{i\theta}|^{4}}u(re^{i\theta})\,d\theta,
$$
which imply
\beq\label{eq-1-lem1}
|\nabla u(0)|
&=&\left[\left |\frac{1}{r\pi}\int_{0}^{2\pi}u(re^{i\theta})\cos\theta \,d\theta
\right |^{2}+\left |\frac{1}{r\pi}\int_{0}^{2\pi}u(re^{i\theta})\sin\theta
\,d\theta \right|^{2}\right]^{1/2}\\ \nonumber
&\leq&\frac{1}{r\pi}\int_{0}^{2\pi}(|\cos\theta|+|\sin\theta|)|u(re^{i\theta})|\,d\theta\\
\nonumber
&\leq&\frac{\sqrt{2}}{r\pi}\int_{0}^{2\pi}|u(re^{i\theta})|\,d\theta.
\eeq
Similar argument shows that
\be\label{eq-x3}
|\nabla v(0)|\leq\frac{\sqrt{2}}{r\pi}\int_{0}^{2\pi}|v(re^{i\theta})|\,d\theta.
\ee
By (\ref{eq-1-lem1}), (\ref{eq-x3}) and Lemma \ref{Lem1}, we obtain that
\begin{eqnarray*}
\Lambda_{f}(0)&\leq&|\nabla u(0)|+|\nabla v(0)|\\
&\leq&\frac{\sqrt{2}}{r\pi}\int_{0}^{2\pi}\big(|u(re^{i\theta})|+|v(re^{i\theta})|\big)\,d\theta\\
&\leq&\frac{2}{r\pi}\int_{0}^{2\pi}|f(re^{i\theta})|\,d\theta.
\end{eqnarray*}
Finally, the desired conclusion follows if we apply the last inequality to the function $F(z)=f(a)-f(z+a)$.
 \epf

\subsection{Proof of Theorem \ref{thm1}} First, we show the ``if" part. By Lemma \ref{lem1}, we have
$$\Lambda_{f}(z)\leq\frac{2}{\pi
\rho}\int_{0}^{2\pi}|f(z)-f(z+\rho e^{i\theta})|\,d\theta,
$$
where $\rho\in(0,d(z)]$, which gives
$$\int_{0}^{r}\Lambda_{f}(z)\rho^{2} \,d\rho\leq\frac{2}{\pi
}\int_{0}^{r}\Big(\rho\int_{0}^{2\pi}|f(z)-f(z+\rho
e^{i\theta})|\,d\theta\Big)\,d\rho,
$$
whence
\begin{eqnarray*}
\Lambda_{f}(z)&\leq&\frac{6}{\pi
r^{3}}\int_{\mathbb{D}(z,r)}|f(z)-f(\zeta)|\,dA(\zeta)\\
&=&\frac{6}{r|\mathbb{D}(z,r)|}\int_{\mathbb{D}(z,r)}|f(z)-f(\zeta)|\,dA(\zeta)\\
&\leq&\frac{6 Mk(r)}{r} =6M\omega\left(\frac{1}{d(z)}\right),
\end{eqnarray*}
where $r=d(z)$.

Next, we prove the ``only if" part. For $z,\
w\in\mathbb{D}$ and $t\in(0,1)$, we have
$$d\big(z+t(w-z)\big)=1-|z+t(w-z)|\geq d(z)-t|w-z|.
$$
If $d(z)-t|w-z|>0$, then
\begin{eqnarray*}
|f(z)-f(w)|&\leq&\left|\int_{0}^{1}\frac{df}{dt}(z+t(w-z))\,dt\right|\\
&\leq&|w-z|\int_{0}^{1}\Lambda_{f}(z+t(w-z))\,dt\\
&\leq&M|w-z|\int_{0}^{1}\omega\left(\frac{1}{d(z)-t|w-z|}\right)dt\\
&=&M\int_{0}^{|w-z|}\omega\left(\frac{1}{d(z)-t}\right)dt.
\end{eqnarray*}
Hence
\begin{eqnarray*}
\frac{1}{|\mathbb{D}(z,r)|}\int_{\mathbb{D}(z,r)}|f(\zeta)-f(z)|\,dA(\zeta)
&\leq&
\frac{M}{|\mathbb{D}_{r}|}\int_{\mathbb{D}_{r}}\left\{\int_{0}^{|\xi|}\omega\left(\frac{1}{d(z)-t}\right)dt\right\}dA(\xi)\\
&=&\frac{2M}{r^{2}}\int_{0}^{r}\rho\left\{\int_{0}^{\rho}\omega\Big(\frac{1}{d(z)-t}\Big)\,dt\right\}d\rho\\
&\leq& \frac{2M}{r^{2}}\int_{0}^{r}\left(\int_{t}^{r}\rho
\,d\rho\right)\omega\left(\frac{1}{r-t}\right)dt\\
&\leq&\frac{2M}{r}\int_{0}^{r}(r-t)\omega\left(\frac{1}{r-t}\right)dt\\
&\leq&\frac{2M}{r}r\omega\Big(\frac{1}{r}\Big)\int_{0}^{r}\,dt\\
&=& 2Mr\omega\big(\frac{1}{r}\big).
\end{eqnarray*}
The proof of this theorem is complete. \qed \medskip

The following result from \cite{KM} is needed in the proof of Theorem \ref{Thm-B}.

\begin{Lem}{\rm \cite[Theorem 1]{KM}}\label{Lem2.0}
Let $u$ be a real harmonic function of $\mathbb{D}$ into $(-1,1)$.
Then for $z\in\mathbb{D}$, the following sharp inequality  holds:
$$|\nabla u(z)|\leq\frac{4}{\pi}\frac{1-u^{2}(z)}{1-|z|^{2}}.
$$
\end{Lem}


\subsection{Proof of Theorem \ref{Thm-B}}
Without loss of generality, we assume that $f$ is not constant. The
implication $f\in L_{\omega}(G)\Rightarrow |f|\in
L_{\omega}(G)\Rightarrow |f|\in L_{\omega}(G,\partial G)$ is
obvious, and so we only need to prove the implication $|f|\in
L_{\omega}(G)\Rightarrow f\in L_{\omega}(G).$
  For a fixed $z\in G$, let
$$M_{z}=\sup\{|f(\zeta)|:~|\zeta-z|<d_{G}(z)\}
$$
and for $\xi\in\mathbb{D}$,
$$T_{f}(\xi)=f(z+d_{G}(z)\xi)/M_{z}.
$$
Obviously, $|T_{f}(\xi)|<1$ and thus Lemma
\Ref{Lem2.0} implies that
$$|\nabla T_{f}(\xi)|\leq\frac{4}{\pi}\left(\frac{1-T_{f}^{2}(\xi)}{1-|\xi|^{2}}\right),
$$
which gives
$$\frac{d_{G}(z)|\nabla f(z)|}{M_{z}}
=|\nabla T_{f}(0)|\leq\frac{4}{\pi}\left(1-\frac{f^{2}(z)}{M^{2}_{z}}\right)\leq\frac{8}{\pi}\left(1-\frac{|f(z)|}{M_{z}}\right),
$$
that is,
\be\label{eq-xx1}
d_{G}(z)|\nabla f(z)|\leq\frac{8}{\pi}\left(M_{z}-|f(z)|\right).
\ee

For a fixed $\varepsilon_{0}>0$, there exists a $\zeta\in\partial G$
such that $|\zeta-z|<(1+\varepsilon_{0})d_G(z)$.  Then, for
$w\in\mathbb{D}(z,d_G(z))$, we have
\begin{eqnarray*}
|f(w)|-|f(z)|&\leq&\big ||f(w)|-|f(\zeta)|\big |+\big ||f(\zeta)|-|f(z)|\big |\\
&\leq&M\omega((2+\varepsilon_{0})d_G(z))+M\omega((1+\varepsilon_{0})d_G(z)).
\end{eqnarray*}
Now we take $\varepsilon_{0}=1$.
It follows that
$$\sup_{w\in\mathbb{D}(z,d_G(z))}(|f(w)|-|f(z)|)\leq M\big(\omega(3d_G(z))+\omega(2d_G(z))\big)\leq5M\omega(d_G(z))
$$
whence
\be\label{eq-xx2}
M_{z}-|f(z)|\leq 5M\omega(d_G(z)).
\ee
By (\ref{eq-xx1}) and (\ref{eq-xx2}), we conclude that
\be\label{eq-xx3}
|\nabla f(z)|\leq \frac{40M}{\pi}\frac{\omega(d_{G}(z))}{d_{G}(z)}.
\ee

Finally, for $z_{1},z_{2}\in G$, by \cite{L}, there must exist a
rectifiable curve $\gamma$ in $G$ which joins $z_{1}$ and $z_{2}$,
and satisfies (\ref{eq1.0}). Integrating (\ref{eq-xx3}) along
$\gamma$, we obtain that
$$|f(z_{1})-f(z_{2})|\leq\int_{\gamma}|\nabla f(\zeta)|\,ds(z)
\leq\frac{40M}{\pi}\int_{\gamma}\frac{\omega(d_{G}(z))}{d_{G}(z)}\,ds(z)\leq
C\omega(|z_{1}-z_{2}|),
$$
where $C$ is a constant. The proof of this theorem is complete.
\qed

\subsection*{Proof of Theorem \ref{Thmy1}}
 The implication $g,\ h\in
L_{\omega}(G)\Longleftrightarrow |g|,\ |h|\in L_{\omega}(G)$ follows
from Theorem \Ref{ThmB}.  We only need to prove $f\in
L_{\omega}(G)\Longrightarrow g,\ h\in L_{\omega}(G),$ because the
implication $g,\ h\in L_{\omega}(G)\Longrightarrow f\in
L_{\omega}(G) $ is obvious. Let $f=h+\overline{g}$ in
 $G$, where $h$ and $g$ are holomorphic in $G$.
 It is easy to know that $f\in L_{\omega}(G)\Longrightarrow \overline{f}\in L_{\omega}(G).$
This implies that $u=\mbox{Re}f_{1}\in L_{\omega}(G)$ and
$v=\mbox{Im}f_{2}\in L_{\omega}(G)$, where $f_{1}=h+g$ and
$f_{2}=h-g$.

We claim that $f_{1},\ f_{2}\in L_{\omega}(G)$. Now we come to prove
this claim. For a fixed $z\in G$, let
$$M_{z}=\sup\{|u(\zeta)|:~|\zeta-z|<d(z)\}~\mbox{and}~T_{u}(\xi)=\frac{u(z+d(z)\xi)}{M_{z}},~ \xi\in\mathbb{D}.$$
Then for any $\xi\in\mathbb{D}$, $|T_{u}(\xi)|<1$ and by Lemma
\Ref{Lem2.0}, we have $$|\nabla
T_{u}(\xi)|\leq\frac{4}{\pi}\left(\frac{1-T_{u}^{2}(\xi)}{1-|\xi|^{2}}\right).$$
This gives $$\frac{d(z)|\nabla u(z)|}{M_{z}}=|\nabla
T_{u}(0)|\leq\frac{4}{\pi}\left(1-\frac{u^{2}(z)}{M^{2}_{z}}\right)\leq\frac{8}{\pi}\left(1-\frac{|u(z)|}{M_{z}}\right),$$
which yields \be\label{eq-xx-1}d(z)|f_{1}'(z)|=d(z)|\nabla
u(z)|\leq\frac{8}{\pi}\left(M_{z}-|u(z)|\right).\ee 

For a fixed $\varepsilon_{0}>0$, there exists a $\zeta\in\partial G$
such that $|\zeta-z|<(1+\varepsilon_{0})d_G(z)$.  Then, for
$w\in\mathbb{D}(z,d_G(z))$, we have
\begin{eqnarray*}
|u(w)|-|u(z)|&\leq&\big ||u(w)|-|u(\zeta)|\big |+\big ||u(\zeta)|-|u(z)|\big |\\
&\leq&M\omega((2+\varepsilon_{0})d_G(z))+M\omega((1+\varepsilon_{0})d_G(z)).
\end{eqnarray*}
Now we take $\varepsilon_{0}=1$. It follows that
$$\sup_{w\in\mathbb{D}(z,d_G(z))}(|u(w)|-|u(z)|)\leq M\big(\omega(3d_G(z))+\omega(2d_G(z))\big)\leq5M\omega(d_G(z))
$$
whence \be\label{eq-xx-5} M_{z}-|u(z)|\leq 5M\omega(d_G(z)). \ee By
(\ref{eq-xx-1}) and (\ref{eq-xx-5}), we conclude that
\be\label{eq-xx-3} | f_{1}'(z)|\leq
\frac{40M}{\pi}\frac{\omega(d_{G}(z))}{d_{G}(z)}. \ee

Finally, for $z_{1},z_{2}\in G$, by \cite{L}, there must exist a
rectifiable curve $\gamma$ in $G$ which joins $z_{1}$ and $z_{2}$,
and satisfies (\ref{eq1.0}). Integrating (\ref{eq-xx-3}) along
$\gamma$, we obtain that
$$|f_{1}(z_{1})-f_{1}(z_{2})|\leq\int_{\gamma}| f_{1}'(\zeta)|\,ds(z)
\leq\frac{40M}{\pi}\int_{\gamma}\frac{\omega(d_{G}(z))}{d_{G}(z)}\,ds(z)\leq
C\omega(|z_{1}-z_{2}|),
$$
where $C$ is a constant. This gives  $f_{1}\in L_{\omega}(G).$ By
 similar arguments, we know that $f_{2}\in
L_{\omega}(G).$ Hence $(f_{1}+f_{2})\in L_{\omega}(G)$ and
$(f_{1}-f_{2})\in L_{\omega}(G)$. Therefore,
$$h=\frac{f_{1}+f_{2}}{2}\in
L_{\omega}(G)~\mbox{and}~g=\frac{f_{1}-f_{2}}{2}\in L_{\omega}(G).$$
The proof of this theorem is completed. \qed

\section{Estimates on $BMO_{2}$}\label{csw-sec3}

Green's theorem (cf. \cite{CPW4,CPW5}) states that if $g\in
C^{2}(\mathbb{D})$, i.e., twice continuously differentiable in
$\ID$, then
\be\label{eq1.2x}
\frac{1}{2\pi}\int_{0}^{2\pi}g(re^{i\theta})\,d\theta=g(0)+
\frac{1}{2}\int_{\mathbb{D}_{r}}\Delta
g(z)\log\frac{r}{|z|}\,d\sigma(z)
\ee
for $r\in (0, 1)$, where $d\sigma$ denotes the normalized area measure in $\mathbb{D}$.


\begin{lem}\label{lem-2}
For $r\in(0,1)$, let
$$M_{p}^{p}(r,f)=\frac{1}{2\pi}\int_{0}^{2\pi}|f(re^{i\theta})|^{p}\,d\theta,
$$
where $f$ is a harmonic mapping in $\mathbb{D}$.
Then for $p\in[2,\infty)$,  $M_{p}^{p}(r,f)$ is a increasing
function on $r$ in $(0,1)$ and
\beq\label{eq-1-lem2}
r\frac{d}{dr}M_{p}^{p}(r,f) &=& p\int_{\mathbb{D}_{r}}\Big [\big
(\frac{p}{2}-1\big )|f(z)|^{p-4}
|f_{z}(z)\overline{f(z)}+f(z)\overline{f_{\overline{z}}(z)}|^{2}\\
\nonumber && \hspace{1.5cm} +|f(z)|^{p-2}|\widehat{\nabla
f}(z)|^{2}\Big]\, d\sigma(z),
\eeq
where $|\widehat{\nabla f}|=(|f_{z}|^{2}+|f_{\overline{z}}|^{2})^{1/2}$.
\end{lem}
\bpf Since $|f|^{p}$ is subharmonic in $\ID$, we see that
$M_{p}^{p}(r,f)$ is an increasing function on $r$ in $(0,1)$, where
$p\in[2,\infty)$. On the other hand, by (\ref{eq1.2x}), we have
\beq\label{eq-2-lem2} \nonumber
r\frac{d}{dr}M_{p}^{p}(r,f)&=&\frac{1}{2}\int_{\mathbb{D}_{r}}\Delta\big
(|f(z)|^{p}\big )\,d\sigma(z)\\
\nonumber &=&
p\int_{\mathbb{D}_{r}}\Big  [\big (\frac{p}{2}-1\big )|f(z)|^{p-4}
\big |f_{z}(z)\overline{f(z)}+f(z)\overline{f_{\overline{z}}(z)}\big
|^{2}\\ \nonumber && \hspace{1.5cm} +|f(z)|^{p-2}|\widehat{\nabla
f}(z)|^{2}\Big]\, d\sigma(z).
\eeq
The proof of this lemma is complete. \epf

\begin{lem}\label{lem-2.0}
For $r\in(0,1)$ and $p\in[2,\infty), $ let
$$I_{p}(r,f)=\left\{\frac{1}{|\mathbb{D}_{r}|}\int_{\mathbb{D}_{r}}|f(z)|^p\,dA(z)\right\}^{1/p},
$$
where $f$ is harmonic in $\mathbb{D}$. Then the function $I_{p}(r,f)$ is increasing on $r$ in $(0,1)$.
\end{lem}
\bpf Since
\be\label{eq-lem-2-1}
\int_{\mathbb{D}_{r}}|f(z)|^{p}dA(z)=2\pi\int_{0}^{r}\rho
M_{p}^{p}(\rho,f)\,d\rho,
\ee
we see that
\be\label{eq-lem-2-2}
\frac{d}{dr}\int_{\mathbb{D}_{r}}|f(z)|^{p}\,dA(z)=2\pi r M_{p}^{p}(r,f).
\ee
By (\ref{eq-lem-2-1}), (\ref{eq-lem-2-2}) and Lemma \ref{lem-2}, we get
\be\label{eq-lem-2-3}
M_{p}^{p}(r,f)- I_{p}^{p}(r,f)=\frac{1}{|\mathbb{D}_{r}|}\int_{0}^{r}\frac{d}{dt}M_{p}^{p}(t,f)|\mathbb{D}_{t}|\,dt\geq0.
\ee
By (\ref{eq-lem-2-1}), (\ref{eq-lem-2-3}) and elementary computations, we conclude that
\beq\label{eq-4-lem-2.0} \nonumber
\frac{d}{dr}I_{p}^{p}(r,f)&=&\frac{|\mathbb{D}_{r}|\frac{d}{dr}\int_{\mathbb{D}_{r}}|f(z)|^{p}\,dA(z)-
\int_{\mathbb{D}_{r}}|f(z)|^{p}\,dA(z)\frac{d}{dr}|\mathbb{D}_{r}|}{|\mathbb{D}_{r}|^{2}}\\
\nonumber &=&\frac{2\pi
r\left[|\mathbb{D}_{r}|M_{p}^{p}(r,f)-\int_{\mathbb{D}_{r}}|f(z)|^{p}\,dA(z)\right]}{|\mathbb{D}_{r}|^{2}}\\
\nonumber &\geq&0.
\eeq
Hence the function $I_{p}(r,f)$ is increasing on $r$ in $(0,1)$. The proof of this lemma is complete.
\epf

\begin{lem}\label{lem2}
For fixed $a\in\mathbb{D}$, let $\phi_{a}(z)=a+(1-|a|)z$ in $\mathbb{D}$.
Then for $p\in[2,\infty)$,
\be\label{eqx-4}
\|f\|_{BMO_{p}}=\sup_{a\in\mathbb{D}}\left\{\frac{1}{|\mathbb{D}|}
\int_{\mathbb{D}}|f(\phi_{a}(z))-f(\phi_{a}(0))|^{p}\,dA(z)\right\}^{1/p},
\ee
where $f$ is harmonic in $\mathbb{D}$.
\end{lem}
\bpf
It is not difficult to see that
\be\label{eqx-7}
\sup_{a\in\mathbb{D}}\left\{\frac{1}{|\mathbb{D}|}
\int_{\mathbb{D}}|f(\phi_{a}(z))-f(\phi_{a}(0))|^{p}\,dA(z)\right\}^{1/p}\leq\|f\|_{BMO_{p}}.
\ee
On the other hand, by elementary calculations and Lemma \ref{lem-2.0}, we have

\vspace{8pt}
$\ds \left\{\frac{1}{|\mathbb{D}(a,r)|}\int_{\mathbb{D}(a,r)}|f(\zeta)-
f(a)|^{p}dA(\zeta)\right\}^{1/p}$

\beq\label{eq-1-lem3}
\nonumber
&\leq& \left\{\frac{1}{|\mathbb{D}(a,1-|a|)|}\int_{\mathbb{D}(a,1-|a|)}|f(\zeta)
- f(a)|^{p}dA(\zeta)\right\}^{1/p}\\ \nonumber
&=&\left\{\frac{1}{|\mathbb{D}|}\int_{\mathbb{D}}|f(\phi_{a}(\zeta))
- f(\phi_{a}(0))|^{p}\,dA(\zeta)\right\}^{1/p},\\ \nonumber
\eeq
where $r\in(0,1-|a|]$. Then
\be\label{eqx-8}
\|f\|_{BMO_{p}}\leq\sup_{a\in\mathbb{D}}\left\{\frac{1}{|\mathbb{D}|}
\int_{\mathbb{D}}|f(\phi_{a}(z))-f(\phi_{a}(0))|^{p}\,dA(z)\right\}^{1/p}.
\ee
Obviously, (\ref{eqx-4}) follows from (\ref{eqx-7}) and (\ref{eqx-8}). \epf

\begin{lem}\label{lem4}
For each fixed $a\in\mathbb{D}$, let
$\phi_{a}(z)=a+(1-|a|)z$ in $\mathbb{D}$. Then
\be\label{eq-x9}
|\phi_{a}'(z)|\leq\frac{1-|\phi_{a}(z)|^{2}}{1-|z|^{2}}.
\ee
\end{lem}
\bpf It is easy to see that $f$ is analytic and for all
$z\in\mathbb{D}$, $|\phi_{a}(z)|\leq1$. Then (\ref{eq-x9}) follows
from the Schwarz-Pick Lemma. \epf

\subsection{Proof of Theorem \ref{thm-2}} We first prove $\|f\|_{BMO_{2}}\leq\beta_{f}$.  For a
fixed $a\in\mathbb{D}$, let
$$F_{a}(\zeta)=f(\phi_{a}(\zeta))
$$
in $\mathbb{D}$, where $\phi_{a}(\zeta)=a+(1-|a|)\zeta$. By
Lemma \ref{lem4}, we have
\beq\label{eq-1-thm}
\nonumber
\sup_{\zeta\in\mathbb{D}}\left\{(1-|\zeta|^{2})\Lambda_{F_{a}}(\zeta)\right\}
&=&\sup_{\zeta\in\mathbb{D}}\left\{(1-|\zeta|^{2})\Lambda_{f}(\phi_{a}(\zeta))|\phi_{a}'(\zeta)|\right\}\\
\nonumber
&\leq&\sup_{\zeta\in\mathbb{D}}\left\{(1-|\phi_{a}(\zeta)|^{2})\Lambda_{f}(\phi_{a}(\zeta))\right\}\\
\nonumber &\leq&\beta_{f}.
\eeq
Then Lemma \ref{lem-2} leads to
\begin{eqnarray*}
\frac{d}{dr}M_{2}^{2}\big(r,F_{a}(re^{i\theta})-F_{a}(0)\big)&=&\frac{2}{r\pi}\int_{\mathbb{D}_{r}}|\widehat{\nabla
F_{a}}(\zeta)|^{2}\,dA(\zeta)\\
&\leq&\frac{2}{r\pi}\int_{\mathbb{D}_{r}} \Lambda_{F_{a}}^{2}(\zeta)\,dA(\zeta)
\end{eqnarray*}
\begin{eqnarray*}
&\leq&\frac{2\beta_{f}^{2}}{r\pi}\int_{\mathbb{D}_{r}}\frac{dA(\zeta)}{(1-|\zeta|^{2})^{2}}\\
&=&\frac{4\beta_{f}^{2}}{r}\int_{0}^{r}\frac{\rho}{(1-\rho^{2})^{2}}\,d\rho\\
&=&2\beta_{f}^{2}\sum_{n=1}^{\infty}r^{2n-1},
\end{eqnarray*}
which gives
$$M_{2}^{2}(r,F_{a}(re^{i\theta})-F_{a}(0))\leq\beta_{f}^{2}\sum_{n=1}^{\infty}\frac{r^{2n}}{n}.
$$
Since
\begin{eqnarray*}
\int_{0}^{1}2rM_{2}^{2}\big(r,F_{a}(re^{i\theta})-F_{a}(0)\big)\,dr
&=&\frac{1}{\pi}\int_{0}^{1}\int_{0}^{2\pi}r|F_{a}(re^{i\theta})-F_{a}(0)|^{2}\,d\theta \,dr\\
&=&\frac{1}{|\mathbb{D}|}\int_{\mathbb{D}}|F_{a}(\zeta)-F_{a}(0)|^{2}dA(\zeta),
\end{eqnarray*}
we see that
\begin{eqnarray*}
\frac{1}{|\mathbb{D}|}\int_{\mathbb{D}}|F_{a}(\zeta)-F_{a}(0)|^{2}dA(\zeta)
\leq\int_{0}^{1}2\beta_{f}^{2}\sum_{n=1}^{\infty}\frac{r^{2n+1}}{n}dr
=\beta_{f}^{2}\sum_{n=1}^{\infty}\frac{1}{n(n+1)}
=\beta_{f}^{2},
\end{eqnarray*}
whence
$$\|f\|_{BMO_{2}}\leq\beta_{f}.
$$

Next, we prove $\beta_{f}\leq2\|f\|_{BMO_{2}}.$ By Lemma \ref{lem-2}
and the subharmonicity of $|\widehat{\nabla F_{a}}|^{2}$, we have
\begin{eqnarray*}
\frac{2}{r}\int_{0}^{r}\rho|\widehat{\nabla
F_{a}}(0)|^{2}d\rho&\leq&\frac{2}{r}\int_{0}^{r}\rho
\left[\frac{1}{2\pi}\int_{0}^{2\pi}|\widehat{\nabla F_{a}}(\rho
e^{i\theta})|^{2}d\theta\right] d\rho\\
&=& \frac{1}{r\pi}\int_{\mathbb{D}_{r}}|\widehat{\nabla
F_{a}}(\zeta)|^{2}dA(\zeta)\\
&=&
\frac{1}{2}\frac{d}{dr}M_{2}^{2}(r,F_{a}(re^{i\theta})-F_{a}(0)),
\end{eqnarray*}
which implies
$$|\widehat{\nabla F_{a}}(0)|^{2}r^{2}\leq M_{2}^{2}(r,F_{a}(re^{i\theta})-F_{a}(0)).
$$ It follows that
\begin{eqnarray*}
\frac{|\widehat{\nabla F_{a}}(0)|^{2}}{4}&=&\int_{0}^{1}|\widehat{\nabla
F_{a}}(0)|^{2}r^{3}dr\leq\frac{1}{2\pi}\int_{\mathbb{D}}|F_{a}(\zeta)-F_{a}(0)|^{2}dA(\zeta),
\end{eqnarray*}
whence
\be\label{eq1h}
\frac{\Lambda_{F_{a}}^{2}(0)}{4}\leq\frac{|\widehat{\nabla
F_{a}}(0)|^{2}}{2}\leq\frac{1}{|\mathbb{D}|}\int_{\mathbb{D}}|F_{a}(\zeta)-F_{a}(0)|^{2}dA(\zeta).
\ee
On the other hand,
\be\label{eq2h}
\beta_{f}\leq\sup_{a\in\mathbb{D}}\Lambda_{F_{a}}(0).
\ee
By (\ref{eq1h}) and (\ref{eq2h}), we have
$$\beta_{f}\leq2\|f\|_{BMO_{2}}.
$$

It remains to prove the sharpness in the inequalities. Obviously, the equality sign in the first inequality of
(\ref{eq-a}) occurs when $f$ is constant. For the sharpness part of the second inequality of (\ref{eq-a}),
we let
$$f(z)=C(z+\overline{z}),
$$
where $C$ is a constant. Then
$$\beta_{f}=\sup_{z\in\mathbb{D}}\{(1-|z|^{2})\Lambda_{f}(z)\}=2|C|
$$
and
\begin{eqnarray*}
\|f\|_{BMO_{2}}&=&\sup_{a\in\mathbb{D}}\left\{\frac{1}{|\mathbb{D}|}
\int_{\mathbb{D}}|F_{a}(z)-F_{a}(0)|^{2}dA(z)\right\}^{\frac{1}{2}}\\
&=&|C|\sup_{a\in\mathbb{D}}\left\{\frac{1}{|\mathbb{D}|}
\int_{\mathbb{D}}(1-|a|)^{2}|z+\overline{z}|^{2}dA(z)\right\}^{\frac{1}{2}}\\
&=&|C|\sup_{a\in\mathbb{D}}\left\{\frac{4(1-|a|)^{2}}{\pi}\int_{0}^{1}\int_{0}^{2\pi}r^{3}\cos^{2}\theta
d\theta dr\right\}^{\frac{1}{2}}\\
&=&|C|\sup_{a\in\mathbb{D}}(1-|a|)\\
&=&|C|,
\end{eqnarray*}
whence
$$\beta_{f}=2\|f\|_{BMO_{2}}.
$$
The proof of this theorem is complete.\qed

\subsection*{\sc Acknowledgments}
The research of Matti Vuorinen was supported by the Academy of
Finland, Project 2600066611. The research of Sh. Chen and X. Wang
was partly supported by NSF of China (No. 11071063). The authors
wish to thank their gratitude to professor Miroslav Pavlovic for
valuable comments.

\end{document}